\documentclass{amsart}
\usepackage{amscd}
\usepackage{amssymb}

\input xy
\xyoption{all}
\newcommand{\g}{\mathfrak{g}}

\newcommand{\gl}{\mathfrak{gl}}
\newcommand{\rip}{\mathfrak{sl}}
\newcommand{\ut}{\mathfrak{u}}
\newcommand{\lt}{\mathfrak{l}}

\newcommand{\iso}{\overset{\sim}{\to}}

\DeclareMathOperator{\pro}{pro}

\DeclareMathOperator*{\coli}{colim}

\DeclareMathOperator{\coker}{coker}

\newcommand{\ev}{ev}

\newcommand{\Z}{\mathbb{Z}}

\newcommand{\C}{\mathcal{C}}

\theoremstyle{plain}
\newtheorem{thm}{Theorem}[section]

\newtheorem{lem}[thm]{Lemma}

\theoremstyle{definition}

\newtheorem{rem}[thm]{Remark}

\newtheorem{exa}[thm]{Example}

\theoremstyle{remark}
\newtheorem{notation}[thm]{Notation}

\begin{document}

\title[Cyclic homology and Lie algebra homology]{Cyclic homology of $H$-unital (pro-)algebras,  Lie algebra homology of matrices and a paper of Hanlon's}

\author[Guillermo Corti\~nas]{Guillermo Corti\~nas*}

\address{Departamento de \'Algebra, Geometr\'\i a y Topolog\'\i a\\
          Universidad de Valladolid}

\email{gcorti@agt.uva.es}


\begin{abstract} We consider algebras over a field $k$ of characteristic zero.
The article is concerned with the isomorphism of graded
vectorspaces
\[
H(\gl(A))\iso\wedge (HC(A)[-1])
\]
between the Lie algebra homology of matrices and the free graded commutative
algebra on the cyclic homology of the $k$-algebra $A$, shifted down one 
degree. For unital
algebras this isomorphism is a classical result obtained by Loday and Quillen
and independently by Tsygan.
For $H$-unital algebras, it is
known to hold too, as is that the proof follows from results of Hanlon's.
However, to our knowledge, the proof is not immediate, and has not been 
published. In this paper we fill this gap in the literature by offering a detailed proof. Moreover we establish the isomorphism in the general setting of 
($H$-unital) pro-algebras. 
\end{abstract}

\thanks{(*) Partially supported by CONICET and the 
Ram\'on y Cajal fellowship and by grants ANPCyT PICT 03-12330 and MTM00958.}

\maketitle

\section{Introduction}
We consider not necessarily unital algebras over a fixed field $k$ of characteristic
zero. If $A$ is an algebra, we write $\gl_n A$ for the Lie algebra of $n\times n$ matrices, put
$\gl A=\bigcup_n\gl_nA$ and consider the Lie algebra homology $H(\gl A)$ 
and the cyclic homology $HC(A)$. We define a natural map of graded vectorspaces
\begin{equation}\label{varfintro}
\varphi: H(\gl A)\to \Lambda (HC(A)[-1])
\end{equation}
Recall $A$ is called $H$-unital if its bar homology vanishes
(\cite{wod}). We prove here that $\varphi$ is an isomorphism if $A$ is $H$-unital.
We obtain this as a particular case of a more general result concerning pro-algebras. In this paper a pro-algebra is an inverse system 

\begin{equation}\label{sigmas}
A=\{A_{n+1}\to A_n\}_n
\end{equation}
of algebras and algebra homomorphisms indexed by the positive integers. Note any algebra $A$ can be regarded as a constant pro-algebra where all the maps \eqref{sigmas} are identity maps.
A pro-algebra $A$ is $H$-unital if for each $r\ge 1$, $H_r^{bar}A=\{H_r^{bar}A_n\}_n$ is the zero pro-vectorspace; this means
that for each $r\ge 1$ and each $n$ there is an $m=m(n,r)\ge n$ such that the
map
\[
H^{bar}_r(A_m)\to H^{bar}_r(A_n)
\]
is zero. The map \eqref{varfintro} induces, in each degree $r$ and each level $n$, a map
$\varphi_{r,n}:H_r(\gl A_n)\to (\Lambda (HC(A_n)[-1]))_r$. 
As $n$ varies, we obtain a map of pro-vectorspaces
\begin{equation}\label{homofi}
\varphi_r:H_r(\gl A)\to (\Lambda (HC(A)[-1]))_r
\end{equation}
We show (in Thm. \ref{hhanwi}) that if $A$ is $H$-unital, then \eqref{homofi} is an isomorphism of pro-vectorspaces; this
means that the pro-vectorspaces $\{\ker\varphi_{r,n}\}_n$ and 
$\{\coker\varphi_{r,n}\}_n$ are zero in the sense explained above. 
In particular, applying this to constant pro-algebras, we get that 
\eqref{varfintro} is an isomorphism for any $H$-unital algebra $A$. 

\smallskip

The isomorphism
\begin{equation}\label{iso}
H(\gl A)\cong \Lambda (HC(A)[-1])
\end{equation}
for a unital algebra $A$ was proved by Loday and Quillen (\cite{lq}) and independently by Tsygan
(see \cite{t} for the announcement and \cite{ft} for the proof). For $H$-unital $A$ it is cited in the literature
(e.g. in \cite[E.10.2.6]{lod} and \cite[pp 88, line 4]{qs}) as following from the results of
Hanlon's paper \cite{han}. Hanlon's paper deals with (Lie, bar and cyclic) cohomology
of finite dimensional $\mathbb{C}$-algebras. What is immediate from the results of \cite{han}
is that if $A$ is finite dimensional and $k=\mathbb{C}$ then there is a map between the dual
spaces
\begin{equation}\label{varpsi}
(\Lambda (HC(A)[-1]))^\vee\to (H(\gl A))^\vee
\end{equation}
which is an isomorphism if $A$ is $H$-unital. The transpose of \eqref{varpsi}
gives an isomorphism $\Lambda (HC(A)[-1])\iso H(\gl A)$ for $A$ finite
dimensional over $\mathbb{C}$. The isomorphism \eqref{iso} for 
$H$-unital algebras of not necessarily finite dimension over any field $k$ of characteristic zero also follows from the material in \cite{han}, but the proof is more delicate. The purpose of this paper is to give a detailed proof of 
\eqref{iso} in the general case of $H$-unital pro-algebras over any field $k$
of characteristic zero.

\smallskip
Next we explain the main ideas of the proof of \eqref{iso} in the case
of $H$-unital algebras. The case of general pro-algebras is proved in a similar manner;
see the proof of Thm \ref{hhanwi} for details. The first step, developed in Section 1, is an elementary application of representation
theory, (mainly of the highest weight theorem \cite[Thm 14.13]{hum2}). The Lie
algebra $\gl_n(k)$ acts naturally on the Chevalley-Eilenberg complex
$C(\gl_nA)$; under this action, the latter complex decomposes into 
isotypic components 
\begin{equation}\label{decintro}
C(\gl_nA)=\bigoplus_\mu C(\gl_n A)_\mu
\end{equation}
Sitting inside each component $C(\gl_n A)_\mu$ is the subcomplex of maximal
weight vectors 
\begin{equation}\label{inclusion}
M_\mu C(\gl_nA)
\subset C(\gl_nA)_\mu
\end{equation}
It turns out that $M_\mu C(\gl_n A)$ and  $H(M_\mu C(\gl_nA))$ generate
$C(\gl_n A)_\mu$ and

\noindent $H(C(\gl_n A)_\mu)$ as $\gl_nk$-modules (see lemmas
\ref{lemuno}, \ref{lem2}). 
Moreover in the particular case of the isotype of the trivial representation,
the inclusion \eqref{inclusion} is an equality (see \ref{lemuno}). 
In Section \ref{apli} we construct a chain homomorphism (see Thm \ref{hanwi})
\begin{equation}\label{fintrito}
\phi: M_{\mu}C(\gl_nA)
\to \Lambda (C^\lambda A[-1])\otimes C^{bar,\mu}A
\end{equation}
where $C^\lambda A$ is Connes' complex for cyclic homology, and $C^{bar,\mu}A$
is a certain complex which depends on $\mu$ and on the bar complex 
of $A$ (see Thm \ref{hanwi} for the full expression of the target of the map 
\eqref{fintrito}). If $\mu$ is the isotype of the trivial representation, $C^{bar,\mu}A$ is just $k[0]$; otherwise it is a complex which is acyclic if $A$ is $H$-unital.
We show (also in Thm \ref{hanwi}) 
that for each $p$ there is an $n_0=n_0(p)$ such that for all $\mu$ and all $n\ge n_0$, the
map $\phi_p$ is an isomorphism
\begin{equation}\label{fipe}
\phi_p:M_{\mu}C_p(\gl_nA)\iso (\Lambda (C^\lambda A[-1])\otimes 
C^{bar,\mu}A)_p
\end{equation}
By what we have explained above, this implies that for $A$ $H$-unital,
and $n\ge n_0(p+1)$, the homology 
of $M_{\mu}C(\gl_nA)$, and therefore also that of 
$C(\gl_n A)_\mu$, vanishes in degree $p$ for all isotypes but that of the trivial representation, in which case it is isomorphic to $(\Lambda (HC(A)[-1]))_p$. 
This proves \eqref{iso} for algebras.

\smallskip

Hanlon's results come up in the proof of \eqref{fipe}.
The proof takes advantage of a certain duality between $\phi$ and a natural
map $\psi$ in the opposite direction defined by Hanlon in \cite[3.6]{han}, 
which is not a chain map, but is such that $\psi_p$ is an isomorphism for each $p$ 
and each $n$ sufficiently large. See the proof of \ref{hanwi} for details.

\smallskip

The rest of this paper is organized as follows. In Section \ref{apli} we 
apply basic representation theory (essentially the highest weight theorem
\cite[Thm 14.13]{hum2}) to prove some elementary facts concerning the 
decomposition of $C(\gl_n A)$ into isotypic components under the action of 
$\gl_nk$. The main result of Section \ref{lamain} is Theorem \ref{hanwi}, where
the chain map \eqref{fintrito} is constructed and the isomorphism \eqref{fipe}
proved.  
In Section \ref{prosec} we prove the main theorem of the paper (Thm \ref{hhanwi}), 
which says that \eqref{homofi} is an isomorphism for every $H$-unital 
pro-algebra $A$.

\section{Application of representation theory}\label{apli}

Throughout this paper $k$ will be a fixed field of characteristic zero; tensor products,
vectorspaces and algebras are over $k$, as are the various (Lie, cyclic, bar) homologies considered.
If $m\ge 0$ then by a {\it partition} of $m$ of length $l$ we understand a nonincreasing sequence
$\alpha_1\ge\dots\ge\alpha_l$ of positive
integers such that $\sum_i\alpha_i=m$. We write $\emptyset$ for empty partition; it is a partition of $0$ of length $0$.
The set of all partitions of a given $m\ge 1$ is denoted $P(m)$. If $\alpha$ and $\beta$ are partitions of $m$ of lengths $l_1$ and $l_2$ and
$l_1+l_2\le n$,
we put
$$
[\alpha,\beta]_n:=(\alpha_1,\dots,\alpha_{l_1},0,\dots,0,-\beta_{l_2},\dots,-\beta_{1})\in\Z^n.
$$
Write $\gl_nk$ for the Lie algebra of $n\times n$ matrices, and
$U(\gl_nk)$ for the universal enveloping algebra. If $V$ is a
$\gl_nk$-module and $S\subset V$ a subset, then $U(\gl_nk)\cdot S$
will denote the $\gl_nk$-submodule of $V$ generated by $S$. For
$\mu\in k^n$, put
\begin{align*}
w_\mu(V):=&\{v\in V: e_{pp}\cdot v=\mu_p v\quad ( 1\le p\le n)\}\\
M_\mu(V):=&\{v\in w_\mu V: e_{ij}\cdot v=0 \quad (1\le i<j \le n)\}\\
V_\mu:=&U(\gl_nk)\cdot M_{\mu}(V)
\end{align*}

If $A$ is any vectorspace and $n\ge 1$ we view
$\gl_nA=\gl_nk\otimes A$ as the tensor product of the adjoint and
the trivial $\gl_nk$-actions; thus
\begin{equation}\label{actua}
g\cdot (h\otimes a)=[g,h]\otimes a.
\end{equation}
If furthemore $A$ is equipped with an associative (not necessarily
unital) algebra structure, then $\gl_nA$ is a Lie algebra and the
induced action on the exterior $k$-algebra $\Lambda\gl_nA$ is
compatible with the Chevalley-Eilenberg boundary map
$\partial:\Lambda^*\gl_nA\to \Lambda^{*-1}\gl_nA$, so that the
Chevalley-Eilenberg complex $C(\gl_nA)=(\Lambda\gl_nA,\partial)$
is a complex in the category of $\gl_nk$-modules. If $\g$ is a Lie
algebra, we write $H(\g)$ for the homology of $C(\g)$.

\begin{lem}\label{lemuno}
\item{i)}The Chevalley-Eilenberg complex splits into a direct sum of subcomplexes
\begin{equation}\label{isodecocheva}
C(\gl_nA)=\bigoplus_{m\ge 0}\bigoplus_{\alpha,\beta\in
P(m)}C(\gl_nA)_{[\alpha,\beta]_n}.
\end{equation}
\item{ii)} $M_{[\emptyset,\emptyset]_n}(C( \gl_nA))=C( \gl_nA)_{[\emptyset,\emptyset]_n}$.
\end{lem}
\begin{proof}

Let the symmetric group $\Sigma_m$ act on the $m$-tensor power
$T^m(\gl_nA)$ by permuting the factors. Consider the idempotent
element $\epsilon_m=(1/m!)\sum_{\sigma\in
\Sigma_m}sg(\sigma)\sigma\in k[\Sigma_m]$. We have
$\Lambda^m\gl_nA=\epsilon_m T^m(\gl_nA)$. Thus, because the
actions of $\Sigma_m$ and $\gl_nk$ on $T^m(\gl_nA)$ commute,
$\epsilon_m(T^m\gl_nA_\mu)=\Lambda^m\gl_nA_\mu=C_m(\gl_nA)_\mu$
for all $\mu\in k^n$. We shall show that
\begin{equation}\label{isodecotenso}
T(\gl_nA)=\bigoplus_m\bigoplus_{\alpha,\beta\in P(m)}T(\gl_nA)_{[\alpha,\beta]_n}.
\end{equation}
Note that, because by definition, the action of $\gl_nk$ on $T^mA$
is trivial,
$$
T(\gl_nA)_\mu=\bigoplus_m (T^m(\gl_nk))_\mu\otimes T^mA.
$$
for all $\mu\in k^n$ . Hence we may assume $A=k$. As both sides of
\eqref{isodecotenso} commute with extension of the scalar field, a
scalar extension--descent argument shows we may further restrict
ourselves to the case when $k$ is algebraicaly closed, or even
more specifically to the case $k=\mathbb{C}$. Hence we can invoke
elementary representation theory of $\gl_nk$ for $k$ algebraically
closed --explained in \cite[Chap. 14]{hum2} for the case
$k=\mathbb{C}$-- by which the decomposition \eqref{isodecotenso}
is valid with $T(\gl_nA)$ replaced by any finite dimensional
representation $V$ of $\gl_nk$ on which the identity matrix acts
trivially, and is the decomposition of $V$ into isotypic
components. The decomposition of part i) of the lemma is proved;
compatibility of the latter with the boundary map follows from the
fact that any homomorphism between nonisomorphic irreducible
representations is zero (Schur's lemma). Next note that in part
ii) the multiplicative structure of $A$ does not play any role.
Moreover as both terms in the identity to prove commute with
scalar extension, and since if $k\subset K$ are fields then every
$K$-vectorspace is the scalar extension of a $k$-vectorspace, if
ii) holds for a particular field of characteristic zero then it
holds for all such fields. Hence we may again assume $k$ is
algebraically closed, or, even more specifically, that
$k=\mathbb{C}$. Now ii) is a consequence of the highest weight
theorem \cite{hum2} applied to the trivial representation of
$\rip_n(\mathbb{C})$.
\end{proof}
\begin{lem}\label{lem2}
$$
H(C(\gl_nA)_{[\alpha,\beta]_n})=U(\gl_nA)\cdot
H(M_{[\alpha,\beta]_n}C(\gl_nA))
$$
\end{lem}
\begin{proof}
We may assume that $k$ is algebraically closed. Choose a graded subspace
\[
V\subset M_{[\alpha,\beta]_n}(\ker\partial_{|C(\gl_nA)})=
\ker\partial\cap M_{[\alpha,\beta]_n}(C(\gl_nA))
\]
such that
$$
V\oplus\partial(M_{[\alpha,\beta]_n}(C(\gl_nA))[+1])=
M_{[\alpha,\beta]_n}(\ker\partial_{|C(\gl_nA)})
$$
By the
highest weight theorem (\cite[Prop. 14.13]{hum2}),
\begin{align*}
\ker\partial_{|C(\gl_nA)_{[\alpha,\beta]_n}}=&U(\gl_nk)\cdot M_{[\alpha,\beta]_n}(\ker\partial_{|C(\gl_nA)_{[\alpha,\beta]_n}})\\
                                                =&U(\gl_nk)\cdot V\oplus U(\gl_nk)\cdot\partial(M_{[\alpha,\beta]_n}(C(\gl_nA)[+1])\\
                                                =&U(\gl_nk)\cdot V\oplus \partial((C(\gl_nA)_{[\alpha,\beta]_n}[+1])
\end{align*}
Thus
\begin{align*}
H(\gl_nA)=&U(\gl_nk)\cdot V\\
          \cong  &U(\gl_nk)\cdot H(M_{[\alpha,\beta]_n}(C(\gl_nA))).
\end{align*}
\end{proof}

In the next lemma and below, we shall write  $V^\vee$ for the
(graded) dual of a (graded) vectorspace $V$. If $f:V\to W$ is a
linear transformation, then ${f}^t:W^\vee\to V^\vee$ is its
transpose.

\begin{lem}\label{dual}
Let $V$ be a finite dimensional vectorspace. Equip $(\Lambda
\gl_nV)^\vee$ with the following $\gl_nk$-module structure
$$
(g\cdot \chi)(x):=\chi(g^t\cdot x).
$$
Then for $m\ge 0$ and $\alpha,\beta\in P(m)$ with $l(\alpha)+l(\beta)\le n$,
the canonical restriction map
\[
(\Lambda \gl_nV)^\vee\to (M_{[\alpha,\beta]_n}\Lambda\gl_nV)^\vee
\]
induces an isomorphism
\begin{equation}\label{res}
M_{[\alpha,\beta]_n}(\Lambda\gl_nV)^\vee \iso
(M_{[\alpha,\beta]_n}\Lambda\gl_nV)^\vee.
\end{equation}
\end{lem}
\begin{proof}
We may assume that $k$ is algebraically closed.
Then (see \cite{hum2})
\begin{equation}\label{weide}
\Lambda \gl_nV=\bigoplus_{\mu\in k^n}w_\mu (\Lambda\gl_nV)
\end{equation}
It follows from this that for all $\mu\in k^n$, the restriction map
$(\Lambda \gl_nV)^\vee\to (w_\mu\Lambda \gl_nV)^\vee$ induces
an isomorphism
\begin{equation}\label{wtowdual}
w_\mu((\Lambda\gl_nV)^\vee)\iso (w_\mu(\Lambda\gl_n V))^\vee
\end{equation}
of which the inverse is induced by the transpose of the projection
of $\Lambda \gl_nV$ onto $w_\mu(\Lambda\gl_nV)$. In particular this holds for $\mu=[\alpha,\beta]_n$. Let $\pi_1$
be the projection of $\Lambda \gl_nV$ onto the summand of
\eqref{weide} corresponding to $\mu=[\alpha,\beta]_n$.
Write $\ut_nk$ and
$\lt_nk$ for the Lie subalgebras $\subset\gl_nk$ of upper triangular and
lower triangular matrices. The map
\eqref{wtowdual} sends $M_{[\alpha,\beta]_n}((\Lambda\gl_nV)^\vee)$
onto
\begin{gather}
\{\chi\in (w_\mu(\Lambda\gl_nV))^\vee: u\cdot \chi=0\ \ \forall u\in\ut_nk\}\nonumber\\
=\{\chi\in (w_\mu(\Lambda\gl_nV))^\vee: \chi(lx)=0\ \ \forall l\in\lt_nk,\ \ x\in \gl_nV\}\label{wdual}
\end{gather}
Consider
the decomposition
\begin{equation}\label{wm}
w_{[\alpha,\beta]_n} (\Lambda\gl_nV)=
\bigoplus_{m\ge 0}
\bigoplus_{\alpha',\beta'\in P(m)}
\Lambda\gl_nV_{[\alpha',\beta']_n}\cap w_{[\alpha,\beta]_n}(\Lambda\gl_nV)
\end{equation}
Write $\pi_2$ for the projection of $w_{[\alpha,\beta]_n}
(\Lambda\gl_nV)$ onto the summand corresponding to
$[\alpha,\beta]_n$. If $(\alpha',\beta')\ne (\alpha,\beta)$ then, by \cite[14.16]{hum2},
every element of $\Lambda\gl_nV_{[\alpha',\beta']_n}\cap
w_{[\alpha,\beta]_n}(\Lambda\gl_nV)$ is a linear combination of
elements of $M_{[\alpha',\beta']_n}\Lambda\gl_nV$ with
coefficients which are
products in $U(\gl_nk)$ of strictly lower triangular matrices. It
follows that $\pi_2^t$ induces an isomorphism between
$(M_{[\alpha,\beta]_n}\Lambda\gl_nV)^\vee$ and \eqref{wdual}.
Hence the map of the lemma is an isomorphism with inverse induced by
$(\pi_2\pi_1)^t$.
\end{proof}

\section{Stable calculation of $M_{[\alpha,\beta]_n}C_p(\gl_nA)$.}\label{lamain}

In preparation for the theorem below, we introduce some notation.
If $V$ is a graded vectorspace, each tensor power $T^mV$ will be considered as
a $\Sigma_m$-module with the action in which the transposition $(i,j)$ acts
as follows
\begin{multline}\label{gradedact}
(i,j)\cdot v_1\otimes\dots\otimes v_i\otimes\dots\otimes v_j\otimes\dots\otimes v_m=\\
(-1)^{|v_i||v_j|}v_1\otimes\dots\otimes v_j\otimes\dots\otimes v_i\otimes\dots\otimes v_m.
\end{multline}
We regard $\Z_m$ embedded
in $\Sigma_m$ through the monomorphism sending the generator $1$ to $(1,\dots,m)$. If $A$ is an algebra, we write
\[
C_n^\lambda A:=(T^{n+1}A)_{\Z/n+1}
\]
and $C^\lambda A=\bigoplus_n C^\lambda_nA$. The {\it Connes cyclic complex} 
(\cite[2.1.4]{lod}) of $A$, $(C^\lambda A,b)$, is a chain
complex whose underlying graded vectorspace is $C^\lambda A$; its homology is the cyclic
homology of $A$, $HC(A)$. Let $T^+A$ be the positive degree part of the tensor algebra. The {\it bar complex} of $A$
is a chain complex $C^{bar}A=(T^+A,b')$ of which the underlying graded
vectorspace is $T^+(A[-1])$ (\cite[Section 2]{wod}). Note that in particular \eqref{gradedact} defines an action of $\Sigma_m$ on $C^{bar}_m(A)$ for each
$m\ge 1$. The homology of $C^{bar}A$ is called the {\it bar homology} of $A$,
denoted $H^{bar}(A)$. Note that, as we do not require our algebras to be unital, any vectorspace $V$ may be regarded
as an algebra in our sense with zero multiplication map; hence $C(\gl_n(V))$, $C^\lambda V$ and $C^{bar}V$ are defined, and are chain
complexes with zero boundary map.

In the next theorem and below, by the canonical inclusion $\gl_nA\subset \gl_{n+1}A$
we mean that given by the monomorphism
\[
g\mapsto \left[\begin{matrix} g&0\\ 0&0\end{matrix}\right].
\]
Let $\alpha$ be a partition of $m\ge 0$. Recall (\cite[pp. 57]{hum2}) that a
{\it standard $\alpha$-tableau} is a filling of the Young diagram of $\alpha$
with the numbers $1,\dots,m$ such that all rows and columns are increasing. We
write $V^\alpha$ for the Specht $k[\Sigma_m]$-module (\cite[4.47]{hum2}) associated to a partition
$\alpha$ of $m$; $V^\alpha$ has a basis with one element
for each standard tableau of shape $\alpha$. In what follows, we shall abuse
notation and identify each standard tableau with the basis vector it
corresponds to. We consider $V^\alpha$ as a graded vectorspace concentrated in
degree zero.

If $V$ is a graded vectorspace, the sum of the coinvariant spaces 
\[
\Lambda V:=\bigoplus_m (T^mV)_{\Sigma_m}
\]
is a graded algebra, the
symmetric algebra of the graded vectorspace $V$; it is commutative in the graded sense. Note that if $W$ is a vectorspace, then $\Lambda W[0]=SW$,
the symmetric algebra and $\Lambda W[-1]=\Lambda W$, the exterior algebra. This admittedly ambiguous notation is
the usual one (see for example \cite{lod}), and should arise no confusion.

If $\phi:V\to W$ is a homomorphism of graded vectorspaces and $p\in \Z$, then
by $\phi_p$ we mean the map $\phi_{|V_p}:V_p\to W_p$.

%

\begin{thm}\label{hanwi}
Let $A$ be a $k$-algebra, $n\ge 1$, $m\ge 0$, and $\alpha$ and $\beta$ partitions of $m$ with
$l(\alpha)+l(\beta)\le n$. There is a natural homomorphism of chain complexes
\begin{equation}\label{elphi}
\phi^n:M_{[\alpha,\beta]_n}(C( \gl_nA))\to \Lambda (C^\lambda A[-1])\otimes 
(T^m(C^{bar}A)\otimes V^\alpha\otimes V^\beta))^{\Sigma_m}
\end{equation}
with the following properties

\smallskip

\item{i)} $\phi^n_p$ is a monomorphism for all $p$ and an isomorphism for
$n\ge p+l(\alpha)+l(\beta)-m$. In particular if $n\ge 2p$ then 
$\phi^n_p$ is an isomorphism.

\smallskip

\item{ii)} The canonical inclusion $\gl_nA\subset \gl_{n+1}A$
sends $M_{[\emptyset,\emptyset]_n}(C( \gl_nA))$ into

\smallskip

$M_{[\emptyset,\emptyset]_{n+1}}(C( \gl_{n+1}A))$ and the following diagram 
commutes
\[
\xymatrix{M_{[\emptyset,\emptyset]_n}(C( \gl_nA))\ar[d]\ar[r]^{\phi^n}& \Lambda (C^\lambda A[-1])\\
M_{[\emptyset,\emptyset]_{n+1}}(C(
\gl_{n+1}A))\ar[ur]_{\phi^{n+1}}&}
\]
\end{thm}

\begin{proof}
We shall define a chain map
\begin{equation*}
\phi':C(\gl_nA)\to R'(A):=\Lambda (C^\lambda A[-1])\otimes T^m(C^{bar}A)\otimes V^\alpha\otimes V^\beta
\end{equation*}
The map of the theorem will be the restriction $\phi:=\phi^n$ of $\phi'$ to $M_{[\alpha,\beta]_n}(\gl_nA)$. We shall show
further below that the image of $\phi'$ really lies
in
\[
R(A):=\Lambda (C^\lambda A[-1])\otimes
(T^m(C^{bar}A)\otimes V^\alpha\otimes V^\beta))^{\Sigma_m}
\]
We will find it convenient further on to view $\Lambda (C^\lambda A[-1])$ as
a trivial $\Sigma_m$-module, so that
\begin{align*}
R'(A)^{\Sigma_m}=&R(A)\\
R'(A)_{\Sigma_m}=&\Lambda (C^\lambda A[-1])\otimes
(T^m(C^{bar}A)\otimes V^\alpha\otimes V^\beta))_{\Sigma_m}
\end{align*}
To define $\phi'$ we take advantage of the $DG$-coalgebra structure
\[
\Delta:C(\gl_nA)\to C(\gl_nA)\otimes C(\gl_nA)=
C(\gl_nk\otimes(A\oplus A))
\]
induced by the diagonal map $A\to A\oplus A$. We shall define two auxiliary maps
$\theta:C(\gl_nA)\to \Lambda (C^\lambda (A)[-1])$ and
$\epsilon:C(\gl_nA)\to T^m(C^{bar}A)\otimes V^\alpha\otimes V^\beta$, and put
\begin{equation}\label{fipri}
\phi'=(\theta\otimes \epsilon)\circ\Delta:C(\gl_nA)\to R'(A)
\end{equation}
To define $\theta$, proceed as follows. 
First define $\theta^1: C(\gl_nA)\to
C^\lambda A[-1]$ by
\[
\theta^1(g_1\land\dots\land g_p)=
\sum_{1\le i_1,\dots,i_p\le n}\sum_{\sigma\in \Sigma_p}\frac{1}{p}sg(\sigma)[(g_{\sigma{1}})_{i_1,i_2}\otimes\dots\otimes
(g_{\sigma(p)})_{i_p,i_1}]
\]
Here $[\quad ]$ denotes the class in $C^\lambda A$; $\theta^1$ is a chain homomorphism because it is the composite of the map \cite[(10.2.3.1)]{lod} and the trace map (\cite[1.2]{lod}),
both of which are chain homomorphisms. Next note that, because $C(\gl_nA)$ is a $DG$-coalgebra and
$\partial(C_1(\gl_nA))=0$, the map $\tilde{\Delta}:C(\gl_nA)\to T(C(\gl_nA))$,
\[
\tilde{\Delta}_p=\sum_{q=0}^p\Delta^{(q)}_p:C_p(\gl_nA)\to T(C(\gl_nA))_p
\]
where $\Delta^{(q)}$ is the $q$-fold comultiplication,
is a chain map. Let $\pi:T(C^\lambda A)\to  \Lambda C^\lambda A$ be the projection. Put
\[
\theta=\pi T(\theta^1)\tilde{\Delta}.
\]
Note $\theta$ is a chain homomorphism, since it is a composite of chain homomorphisms.

In preparation for the definition $\epsilon$ we introduce some notation. Write $M_nA$ for $\gl_nA$ considered
as an associative algebra. If $1\le i,j\le n$ and $g=g_1\otimes\dots \otimes g_p\in C^{bar}_p(M_nA)$, put
\begin{align*}
(g)_{ij}:=&\sum_{1\le l_2,\dots,l_p\le n}(g_1)_{i,l_2}\otimes (g)_{l_2,l_3}\otimes\dots\otimes (g)_{l_p,j}\\
\epsilon_{ij}(g):=&(\sum_{\sigma\in \Sigma_p}sg(\sigma)g_{\sigma(1)}\otimes\dots\otimes g_{\sigma(n)})_{ij}
\end{align*}
One checks that $g\mapsto (g)_{ij}$ is a chain homomorphism $C^{bar}(M_nA)\to C^{bar}A$. As $\epsilon_{ij}$ is the
composite of the latter with the antisymmetrization map of \cite[1.3.4]{lod}, it follows that it is a chain homomorphism.
Next if $\gamma$ is a partition
of $m$ of length $\le n$, $z$ a standard tableau of shape $\gamma$, and $1\le i\le m$, then we put $\rho_i(z)$ for the row
of $z$ containing $i$. Thus if $\alpha$ and $\beta$ are partitions of $m$ of length $\le n$ and $x,y$ are standard tableaux
of shapes $\alpha$ and $\beta$, then the following is a chain map $C(\gl_nA)\to T^m(C^{bar}A)$:
\[
\epsilon_{\rho(x),\rho(y)}=(\epsilon_{\rho_1(x),n+1-\rho_1(y)}\otimes\dots\otimes \epsilon_{\rho_m(x),n+1-\rho_m(y)})\circ\Delta^{(m)}
\]
We define
\[
\epsilon:=\sum_{x,y}\epsilon_{\rho(x),\rho(y)}\otimes x\otimes y
\]
Here the sum runs over all pairs $(x,y)$ of a standard
$\alpha$-tableau $x$ and a standard $\beta$-tableau $y$. It is
clear from the definitions just given that \eqref{fipri} is a
chain homomorphism and that in the case $\alpha=\beta=\emptyset$
it is compatible with the inclusions $\gl_nA\subset\gl_{n+1}A$.
Thus ii) is proved. We remark that the product structure of $A$
does not play any role in the formulas defining $\phi'$, and that
it is equivalent to prove the remaining assertions for $A$ or for
the underlying vectorspace of $A$ equipped with the zero
multiplication map. Hence in what follows $A$ will be any
vectorspace, which we shall implicitly regard as an algebra with
trivial multiplication. As any vectorspace is a filtering colimit
of
 finite dimensional ones, and all functors involved preserve such filtering colimits, we may further
assume that $A$ is finite dimensional. According to \cite[Thms 3.6]{han} there is natural map
\begin{equation}\label{psi}
\psi:(R'(A))_{\Sigma_m}\longrightarrow M_{[\alpha,\beta]_n}C(\gl_nA)
\end{equation}
such that $\psi_p$ is surjective for all $p$ and an isomorphism for
$p\ge l(\alpha)+l(\beta)-m$.
We remark that the bound $p\ge l(\alpha)+l(\beta)$ stated in loc. cit. is a
typos; the proof is done for $p\ge l(\alpha)+l(\beta)-m$
(see \cite[page $219$, line $4$]{han}). Hanlon states his theorem for
$A$ an algebra over $k=\mathbb{C}$. However $\psi$ (whose construction is recalled below) is
defined and natural for algebras over any field $k$, and in particular for vectorspaces $A$ considered
as $k$-algebras with zero multiplication map. Moreover, as both its source and its target commute
with extension of the scalar field, and since if $k\subset K$ are fields then every $K$-vectorspace is
the scalar extension of
a $k$-vectorspace, the fact that, for a particular value of $p$, $\psi_p$ is surjective or an isomorphism for vectorspaces over $\mathbb{C}$ implies
it is one for vectorspaces over any field $k$ of characteristic zero.
Write $\psi'$ for the composite of $\psi$ with the projection $R'(A)\to (R'(A))_{\Sigma_m}$. We shall see that the properties of $\phi$ which remain to be proved follow from Hanlon's result
\cite[Thms 3.6]{han} and from a
certain duality between
$\phi'$ and $\psi'$ which we shall establish. For this we need the explicit description of $\psi'$, which we recall next.
The map $\psi':R'(A)\to C(\gl_nA)$ is a composite
of the form
\[
\psi'=\mu\circ (\hat{\theta}\otimes\hat{\epsilon}):R'(A)\to C(\gl_nA)
\]
where \[
\mu:C(\gl_nA)^{\otimes 2}=C(\gl_n\otimes(A\oplus A))\to C(\gl_nA)
\]
is the  multiplication induced by the sum map $A\oplus A\to A$, and
\[\hat{\theta}:\Lambda (C^\lambda A[-1])\to C(\gl_nA)\]
and
\[\hat{\epsilon}:T^m(C^{bar}A)\otimes V^\alpha\otimes V^\beta\to C(\gl_nA)\]
shall be defined presently. First we introduce notation. If $1\le i,j\le n$ and $a_1,\dots, a_p\in A^{\otimes p}$, set
\[
e_{ij}(a_1\otimes\dots\otimes a_p):=\sum_{1\le l_2,\dots,l_p\le n} e_{i,l_2}(a_1)\land\dots\land e_{l_p,j}(a_p)
\]
Let $\hat{\theta}$ be the graded algebra map determined by
\[
\hat{\theta}([a_1\otimes\dots\otimes a_p])=\sum_{1\le l\le n} e_{l,l}(a_1\otimes\dots\otimes a_p)
\]
If $x$ and $y$ are standard tableaux of shapes $\alpha$ and $\beta$, and $c_1\otimes\dots\otimes c_m\in T^m(C^{bar}A)$, put
\[
\hat{\epsilon}(c_1\otimes\dots\otimes c_m\otimes x\otimes y):=e_{\rho_1(x),n+1-\rho_1(y)}(c_1)\land\dots\land e_{\rho_m(x),n+1-\rho_m(y)}(c_m)
\]
Thus $\psi'$ is well-defined and natural. To express the relation between $\psi'$ and $\phi'$ we need more notation.
Consider the bilinear form
\[
<,>:\gl_nk^{\otimes 2}\to k,\qquad <g,h>:=Tr(g^t\cdot h)
\]

Let
$ev:A\to A^{\vee\vee}$ be the canonical evaluation map; as we are assuming $A$ finite dimensional, it is an isomorphism.
Define a vectorspace isomomorphism
\begin{equation}\label{nu}
\nu:\gl_nA\to \gl_n(A^\vee)^\vee=\gl_nk^\vee\otimes
A^{\vee\vee},\quad \nu(g\otimes a)=<g,\quad>\otimes\ev(a)
\end{equation}

Applying the functor $\Lambda$ to $\nu$ and composing with $\Lambda(\gl_nA^\vee)^\vee\to (\Lambda(\gl_nA^\vee))^\vee$ we
obtain an isomorphism $C(\gl_nA)\iso (C(\gl_n(A^\vee)))^\vee$ which we shall still call $\nu$.
Similarly $\ev$ induces
an isomorphism $R'(A)\iso R'(A^\vee)^\vee$ which we call $\ev$ as well.
The map $N:C^\lambda A\to C^\lambda A$,
\[
N([a_1\otimes\dots\otimes a_n])=n[a_1\otimes\dots\otimes a_n]
\]
extends to an algebra isomomorphism $\Lambda C^\lambda A[-1]
\iso \Lambda C^\lambda A[-1]$ which we shall also call $N$.
I claim that the following diagram
commutes
\begin{equation}\label{conmu}
\xymatrix{C(\gl_nA)\ar[d]_{\nu}\ar[rr]^{\phi'}& &R'(A)\ar[d]^\ev\\
          (C(\gl_n(A^\vee)))^\vee\ar[rr]_{(\psi'\circ (N\otimes 1))^t}& &(R'(A^\vee))^\vee}
\end{equation}
Here $\psi'=\psi'_{A^\vee}$ is the homomorphism the natural transformation
$\psi'$ assigns to the vectorspace (or algebra with trivial multiplication)
$A^\vee$.
Assume the claim is true. Then, because $\psi'$ descends to
$\Sigma_m$-coinvariants, it follows that the image of $\phi'$
is contained in the $\Sigma_m$-invariants. One checks that
$\nu$ satisfies the following identity
\begin{equation}\label{nuequi}
\nu([g,h])(r)=\nu(h)([g^t,r])\qquad (g\in \gl_nk, h\in
C(\gl_n A), r\in C(\gl _n(A^{\vee})))
\end{equation}
It follows from \ref{dual} that the composite of $\nu$ with the restriction map
$C(\gl_n(A^\vee)))^\vee \to (M_{[\alpha,\beta]_n}C(\gl_n(A^\vee)))^\vee$ sends
$M_{[\alpha,\beta]_n}C(\gl_nA)$ isomorphically onto

$(M_{[\alpha,\beta]_n}C(\gl_n(A^\vee)))^\vee$.
Since on the other hand, by \cite[Thm 3.6]{han}, $\psi'(R'(A))=M_{[\alpha,\beta]_n}C(\gl_nA)$,
and the induced map $\psi_p: R'_p(A)_{\Sigma_m}\to M_{[\alpha,\beta]_n}C_p(\gl_nA)$ is injective for
$n\ge p+l(\alpha)+l(\beta)-m$, it follows that the restriction $\phi$ of $\phi'$ to $M_{[\alpha,\beta]_n}C(\gl_nA)$ is an injection, and that
$\phi_p$ is surjective for $n\ge p+l(\alpha)+l(\beta)-m$. It only remains to prove the claim that
\eqref{conmu} commutes.
To see this one checks first that the following two diagrams commute
\begin{equation}\label{simples}
\xymatrix{C(\gl_nA)\ar[d]_{\nu}\ar[r]^{\theta^1}& C^\lambda A[-1]\ar[d]^\ev\\
          (C(\gl_n(A^\vee)))^\vee\ar[r]_{(\hat{\theta})^tN}&
(C^\lambda(A^\vee)[-1])^\vee}
\xymatrix{C(\gl_nA)\ar[d]_{\nu}\ar[r]^{\epsilon_{ij}}& C^{bar}A\ar[d]^\ev\\
          (C(\gl_n(A^\vee)))^\vee\ar[r]_{(\hat{\epsilon}_{ij})^t}& (C^{bar}(A^\vee))^\vee}
\end{equation}
Second, one checks that if $\widetilde{\nu\otimes\nu}$ is the composite of $\nu\otimes\nu$ followed by
the natural isomorphism
\[C(\gl_n(A^\vee))^\vee\otimes C(\gl_n(A^\vee))^\vee\
\iso (C(\gl_n(A^\vee))\otimes C(\gl_n(A^\vee)))^\vee\]
then
\begin{equation}\label{straight}
\xymatrix{C(\gl_nA)\ar[d]_\nu\ar[r]^(0.4){\Delta} &C(\gl_nA)\otimes C(\gl_nA)\ar[d]^{\widetilde{\nu\otimes\nu}}\\
          C(\gl_n(A^\vee))^\vee\ar[r]^(0.4){\mu^t} &(C(\gl_n(A^\vee))\otimes C(\gl_n(A^\vee)))^\vee}
\end{equation}
commutes. Because $\phi'$ is determined by $\theta^1$,
the $\epsilon_{ij}$, and the coproduct structure, and $\psi'$ by
$\hat{\theta}_{|C^\lambda A}$, the $\hat{\epsilon}_{ij}$, and the product
structure, the commutativity of \eqref{conmu} follows from that of
\eqref{simples} and \eqref{straight}.
\end{proof}

\begin{rem} The inclusion $\gl_nA\subset \gl_{n+1}A$ sends
\[
M_{[\alpha,\beta]_n}C(\gl_nA)\to w_{[\alpha,\beta]_n}C(\gl_{n+1}A)
\]
As 
\[
w_{[\alpha,\beta]_n}C(\gl_{n+1}A)\not\supset M_{[\alpha,\beta]_{n+1}}
C(\gl_{n+1}A),
\]
there is no analogue of part ii) of \ref{hanwi} for 
$(\alpha,\beta)\ne (\emptyset,\emptyset)$.
\end{rem}

\begin{rem}\label{noestadefi}
The map $\psi=\psi_A$ of \eqref{psi} is defined for every algebra $A$, finite dimensional or not. It is not a chain
homomorphism. For example if $(\alpha,\beta)=(\emptyset,\emptyset)$ and $n=1$, then, under the identification
$\gl_1A=A$, we have
\begin{align*}
\psi(b[a_0\otimes a_1\otimes a_2])=&\psi(a_0a_1\otimes a_2-a_0\otimes a_1a_2+a_2a_0\otimes a_1)\\
                                  =&(a_0a_1)\land a_2-a_0\land a_1a_2+a_2a_0\land a_1\\
                                  \ne&[a_0,a_1]\land a_2-[a_0,a_2]\land a_1+[a_1,a_2]\land a_0\\
                                  =&\partial(a_0\land a_1\land a_2)=\partial\psi(a_0\otimes a_1\otimes a_2)
\end{align*}
On the other hand, if $A$ is finite dimensional, then a choice of basis for $A$
naturally gives rise to isomorphisms of graded vectorspaces
\begin{equation}\label{nolonger}
C(\gl_n A)\iso (C(\gl_nA))^\vee,\ \ C^\lambda A\iso (C^\lambda A)^\vee\text{ and }C^{bar}A\iso (C^{bar}A)^\vee.
\end{equation}
Using these isomorphisms one can view $C(\gl_nA)$, $C^\lambda A$, $C^{bar}A$ and $R'(A)_{\Sigma_m}$ as cochain complexes;
Hanlon shows \cite[3.11]{han} that $\psi$ is a cochain homomorphism in this sense. In the infinite dimensional case we
no longer have the isomorphisms \eqref{nolonger}. However the map $\phi:C(\gl_nA)\to R(A)$ of \ref{hanwi} is always
a chain homomorphism (independently of the dimension of $A$) and therefore its transpose
$\phi^t:R(A)^\vee\to C(\gl_nA)^\vee$ is always a cochain homomorphism.
\end{rem}

\section{Main theorem}\label{prosec}
\begin{notation}\label{nota}
If $\C$ is a category, we write pro-$\C$ for the category of all
inverse systems
$$
\{\sigma_{n+1}:C_{n+1}\to C_n:n\in\Z_{\ge 1}\}
$$
of objects of $\C$. The set of homomorphisms of two pro-objects
$C$ and $D$ in pro$-\C$ is by definition
$$
\hom_{\pro-\C}(C,D)=\lim_n\coli_m\hom_\C(C_m,D_n)
$$
Thus a map $C\to D\in\pro-\C$ is an equivalence class of maps of
inverse systems
\begin{equation}\label{prorepre}
f:\{C_{m(n)}\}_n\to \{D_n\}_n
\end{equation}
where $m:\Z_{\ge 1}\to \Z_{\ge 1}$ satisfies $m(n)\ge n$ for all
$n$. In the particular case when $m$ is the identity we say that
$f$ (and also its equivalence class) is a {\it level map}.
We remark that mapping each object $C\in \C$ to the constant pro-object gives a fully faithful functor
$\C\to pro-\C$.
Because of this, we shall identify the objects of $\C$ with the constant pro-objects.
Assume $\C$ has a zero object. Then we say that a pro-object $X$ is zero if it is isomorphic to
$0$ in pro-$\C$. This means that
\begin{equation}\label{pronul}
\forall n\ \ \exists m\ge n \text{ such that } \sigma_{n+1}\circ\dots\circ\sigma_m=0.
\end{equation}
If $\C$ happens to be abelian, pro-$\C$ is abelian as well, and if \eqref{prorepre} represents
a homomorphism $[f]:C\to D\in \pro-\C$ then $\{\ker C_{m(n)}\to D_n\}_n$ is a kernel and
$\{\coker C_{m(n)}\to D_n\}_n$ a cokernel for $[f]$. In particular $[f]$ is an isomorphism
if and only if both $\{\ker C_{m(n)}\to D_n\}_n$ and $\{\coker C_{m(n)}\to D_n\}_n$ are zero in the sense of \eqref{pronul}.
If $F:\C\to\mathcal{D}$ is a functor and $C=\{C_n\}_n\in \pro-\C$, we write $F(C)$ for the pro-object
$\{F(C_n)_n\}$. In what follows we shall consider pro-objects in the categories of $k$-algebras and
vectorspaces. We call a pro-algebra $A=\{A_n\}_n$ {\it $H$-unital} if for each $r$ the
pro-vectorspace $H_r^{bar}(A)=0$. According to \eqref{pronul},
this means that for all $r,n\ge 1$ there is an $m=m(n,r)$ such that the map
\[
H^{bar}_r(A_{m})\to H^{bar}_r(A_n)
\]
is zero.

\smallskip

If $A$ is a $k$ algebra, we write $\gl A$ for the Lie algebra of all matrices
of finite size
\[
\gl A=\bigcup_{n\ge 1}\gl_nA.
\]

\end{notation}

\begin{thm}\label{hhanwi}
Let $A=\{A_n\}_n$ be a pro-$k$-algebra. Assume that $A$ is $H$-unital.
Then for each $r\ge 1$ there is a natural isomorphism of
pro-vectorspaces
\begin{equation}\label{aproba}
\varphi_r:H_r(\gl A)\iso (\Lambda (HC(A)[-1]))_r.
\end{equation}
\end{thm}
\begin{proof}
Consider the composite
\[
\varphi^n_r:H_r(\gl_nA_l)\to H_r(M_{[\emptyset,\emptyset]}C(\gl_nA_l))
\to (\Lambda (HC(A_l)[-1]))_r.
\]
By by part ii) of \ref{hanwi}, the $\varphi_r^n$ pass to the limit with respect to $n$, so that we do
have a levelwise map from the left to the right hand side
of \eqref{aproba}. We shall show that $\varphi^n_r$ is an isomorphism for
$n\ge 1+2r$. 
Note that if $\alpha,\beta\in P(m)$ and $m>r$, then by part i) of \ref{hanwi},
we have $M_{[\alpha,\beta]_n}C_r(\gl_n(A_l))=0$ for all $l$.
Thus by \ref{lemuno} and \ref{lem2} we have a levelwise
decomposition
\begin{multline}\label{hanpri}
H_r(\gl_nA)=H_r(M_{[\emptyset,\emptyset]_n}(C(\gl_nA)))\oplus\\
\bigoplus_{1\le m\le r}\bigoplus_{\alpha,\beta\in
P(m)}U(\gl_nk)\cdot H_r(M_{[\alpha,\beta]_n} (C(\gl_nA)))
\end{multline}
If $n>r$, then by \ref{hanwi}, $\varphi^n_r$ is onto and its kernel is the term in the second line of \eqref{hanpri}. If further $n\ge 2r+1$, it follows from
\ref{hanwi} and our hypothesis that $H^{bar}_p(A)=0$ for all $p$, that each of the finitely many
summands in the second line of \eqref{hanpri} is the zero pro-vectorspace. This proves that $\ker\varphi^n_r=0$, whence $\varphi_r^n$ is an isomorphism as wanted.
\end{proof}

\begin{exa} Let $A$ be an algebra; the inclusions $\{A^{n+1}\subset A^n\}_n$ define a pro-algebra $A^\infty:=\{A^n\}_n$.
If $A$ is isomorphic to an ideal of a tensor algebra $TV$ (or more generally of any quasi-free algebra) then 
$A^\infty$ is $H$-unital (see \cite[Section 4]{cor}). 
Hence theorem \ref{hhanwi} applies
in this case.
\end{exa}

\end{document}